\documentclass{article}
\usepackage[utf8]{inputenc} 
\usepackage[T1]{fontenc}    

\usepackage{booktabs}       
\usepackage{microtype}      
\usepackage[noend]{algpseudocode}
\usepackage{algorithm}
\usepackage{xspace}

\usepackage{mathtools}
\usepackage{amsthm}
\usepackage{amssymb}

\usepackage{enumitem}
\usepackage{tikz}
\usetikzlibrary{arrows,decorations.pathreplacing,decorations.markings}

\usepackage[square,comma, sort&compress]{natbib}
\usepackage{jmlr2e}
\usepackage{charter}
\usepackage[vvarbb]{newtxmath}

\usepackage[colorlinks,urlcolor=blue,citecolor=blue,linkcolor=blue]{hyperref}
\hypersetup{hypertexnames=false} 
\usepackage[capitalize]{cleveref}[0.21]

\providecommand{\To}{\textbf{to} }

\newcommand{\lmo}{LMO\xspace}

\DeclareMathOperator*{\conv}{conv}
\DeclareMathOperator{\dist}{dist}

\DeclarePairedDelimiterXPP{\innp}[2]{\mathinner\bgroup}{\langle}{\rangle}%
{\egroup}{#1, #2}

\newcommand{\norm}[1]{\left\| #1 \right\|}

\newcommand{\RR}{\mathbb{R}}

\theoremstyle{plain} \numberwithin{equation}{section}
\newtheorem{theorem}{Theorem}[section]
\numberwithin{theorem}{section}
\newtheorem{corollary}[theorem]{Corollary}

\newtheorem{proposition}[theorem]{Proposition}

\theoremstyle{definition}

\newtheorem{remark}[theorem]{Remark}

\newtheorem{observation}[theorem]{Observation}

\DeclareMathOperator*{\argmin}{argmin}
\DeclareMathOperator*{\argmax}{argmax}

\DeclareMathOperator{\vertex}{vert}

\tikzset{
  on each segment/.style={
    decorate,
    decoration={
      show path construction,
      moveto code={},
      lineto code={
        \path [{#1}]
        (\tikzinputsegmentfirst) -- (\tikzinputsegmentlast);
      },
      curveto code={
        \path [{#1}] (\tikzinputsegmentfirst)
        .. controls
        (\tikzinputsegmentsupporta) and (\tikzinputsegmentsupportb)
        ..
        (\tikzinputsegmentlast);
      },
      closepath code={
        \path [{#1}]
        (\tikzinputsegmentfirst) -- (\tikzinputsegmentlast);
      },
    },
  },
  mid arrow/.style={postaction={decorate,decoration={
        markings,
        mark=at position .5 with {\arrow[{#1}]{>}}
      }}}}

\newcommand{\hrulealg}[0]{\vspace{1mm} \hrule \vspace{1mm}}

\title{Alternating Linear Minimization: \\
  inspired by von Neumann's alternating projections}

\author{\name Gábor Braun \email \href{mailto:braun@zib.de}{braun@zib.de} \\
       \addr Zuse Institute Berlin\\
       Berlin, Germany
       \AND
       \name Sebastian Pokutta \email \href{mailto:pokutta@zib.de}{pokutta@zib.de} \\
       \addr
       Zuse Institute Berlin and \\ Institute of Mathematics,
       Technische Universität Berlin\\
       Berlin, Germany
	\AND
	\name Robert Weismantel\email  \href{mailto:robert.weismantel@ifor.math.ethz.ch}{robert.weismantel@ifor.math.ethz.ch}\\
	\addr 
Department of Mathematics, ETH Zürich \\
Zürich, Switzerland
}

\begin{document}

\maketitle

\begin{abstract} 
  In 1933 von Neumann proved a beautiful result that one can approximate a
  point in the intersection of two convex sets by
  alternating projections, i.e., successively
  projecting on one set and then the other. This algorithm assumes
  that one has access to projection operators for both sets.  In this work, we
  consider the much weaker setup where we have only access to linear
  minimization oracles over the convex sets and present an algorithm
  to find a point in the intersection of two convex sets. In the case of polytopes, we obtain an adaptive variant that is exact and decides in finite time whether the two sets intersect and either provides a point in the intersection or a separating hyperplane.
\end{abstract}

\section{Introduction}

We are interested in the following problem and its variants: Given two (compact) convex sets $P$ and $Q$, compute
\begin{equation}
	\label{eq:intersect}
x \in P \cap Q,
\end{equation}
if such an $x$ exists otherwise certify disjointness.  Von
Neumman's \emph{alternating projection algorithm} is a very simple
algorithm to compute such a point in the intersection of $P$ and $Q$:
we alternatingly project onto each of the sets and we are guaranteed
to converge to a point in the intersection if such a point
exists (see \cite{von1949rings,ginat2018method}).
This idea has sparked a lot of follow-up work and eventually
led to splitting algorithms (see e.g., \cite{combettes2011proximal}
for a great overview), where feasible regions are treated separately
and then the results are combined into a convergent scheme. All
these approaches rely on projections onto the convex sets, or
generalizations of projection, so-called proximal operators.

A different more recent line of work considers the composite
optimization framework, where the goal is to
minimize the sum of a smooth function and a non-smooth function.
To solve Problem~\eqref{eq:intersect}, one minimizes the sum of the
convex function \(\norm{x - y}^{2}\) and
the indicator function of \(x \in P\)
(which has value \(0\) if \(x\) is in \(P\) and has value \(+\infty\)
otherwise)
and the indicator function of \(y \in Q\).
Composite optimization is usually solved by smoothening the non-smooth
summands, e.g.,
\citet{argyriou14}
\citep[later rediscovered in][]{yurtsever18} smoothens via the
\emph{Moreau envelope} \citep{moreau65,bc17}.   This smoothed problem
can then be solved with various algorithms, in particular Frank–Wolfe
algorithms as done in \citet{argyriou14}, however computing
gradients of the smoothened function involves a proximal operator.

Most closely related to our work are \citet{willner1968distance} that
compute the distance between two polytopes via the vanilla Frank–Wolfe
algorithm with an explicit line search and \citet{wolfe1976finding}
that computes the nearest point in a polytope using Wolfe's algorithm
which can be understood as a quadratic extension of the original
Frank–Wolfe approach and which is closely related to the minimum norm
point algorithm of \citet{fujishige1990dual}.
Frank–Wolfe algorithms (also called conditional gradients)
have the appeal that they access the feasible region only through
linear minimization,
we refer the interested reader to the survey of \cite{CGFWSurvey2022} for a broad exposition and background.
For our purposes, the Block-coordinate Frank–Wolfe algorithm from
\cite{lacoste2013block} is the starting point, which minimizes a
convex function over a product of convex sets, randomly alternating
among the convex sets (the blocks).
Below we will use its version,
the \emph{Cyclic Block-Coordinate Conditional Gradient algorithm
(CBCG)} of \citet{beck2015cyclic}, which alternates between the
convex sets in a cyclic fashion,
and thus is a direct analogue of alternating projections
with good performance in actual computations \citep{beck2015cyclic}.

While all these lines of work are highly interesting and relevant,
here we are concerned with a much simpler and more elementary
question. What happens if we change the access model to the feasible
regions, i.e., rather than having access to projection operators, we
assume that we only have access to \emph{linear minimization oracles
  (LMO)} for each of the two sets. These oracles, given a linear
objective $c \in \RR^n$ as input, return
\[
\argmin_{x \in P} \innp{c}{x} \qquad \text{and} \qquad \argmin_{y \in Q} \innp{c}{y},
\]
respectively.
Compared to projection operators, linear minimization
tends to be considerably cheaper, which is one of our motivations; see e.g., \citet{CP2021} for a discussion. The
other question is whether we can still compute
a point in the intersection, having only access to $P$ and $Q$ through
their linear minimization oracles.

We answer this question in the affirmative and obtain an
\emph{alternating linear minimization} algorithm, that computes $x \in
P \cap Q$ or certifies disjointness for two such compact convex sets
$P$ and $Q$ whose iteration complexity (up to constants depending on
$P$ and $Q$) is basically identical to von Neumann's algorithm. We do
only one LMO call per alternation rather than approximating the
projection via the Frank–Wolfe algorithm e.g., via a Douglas–Ratchford
scheme as done in \cite{millan2021approximate}, typically requiring
many LMO calls per each approximation of the projection operation.

\subsubsection*{Contribution} Our contribution can be roughly summarized as follows:

\paragraph{Alternating Linear Minimization Algorithm}

We specialize the Cyclic Block-Coordinate Conditional Gradient
algorithm to von Neumann's setting, similar to the specialization of
the vanilla Frank–Wolfe algorithm
in \citet{willner1968distance}, arriving at our \emph{Alternating
  Linear Minimizations (ALM)} (Algorithm~\ref{alg:dist_LMO}). This algorithm
is the direct analog to von Neumann's alternating projection
algorithm,
however replacing each projection with a step of the Frank–Wolfe
algorithm, an update based on a single linear minimization
only. Note, that this is in sharp contrast to approaches that approximately solve the projection problem by means of conditional gradients as done e.g., in \citet{millan2021approximate} or when employing conditional gradient sliding techniques \citep{lan2016conditional}; see also \citet{CGFWSurvey2022} for an overview. Surprisingly, we obtain (up to constants) the same rate of
convergence of $O(1/t)$ for the distance to the intersection
after \(t\) linear minimizations as von
Neumann's original algorithm albeit with a different dependence on the
geometric parameters of the two sets (as expected): we obtain roughly
$O\left(\frac{(D_P + D_Q)^2}{t}\right)$ vs  $O\left(\frac{\dist(y_0, P
    \cap Q)^{2}}{t}\right)$,
where \(D_{p}\) and \(D_{Q}\) are diameters of the sets \(P\) and
\(Q\) and \(y_{0}\) is the starting point of the algorithm chosen
arbitrarily; see Proposition~\ref{prop:setIntersection} for details. 
Moreover, our algorithm certifies non-existence of $x \in P \cap Q$ in
case $P \cap Q = \emptyset$ and generates a certifying hyperplane as detailed in Section~\ref{sec:empty-intersection}. 

\paragraph{Approximate and Exact Intersection for Polytopes}

While the above might appear as a straightforward specialization of
more general results, the implied complexities of deciding
(non-)intersection solely via access to linear minimization oracles
(which take a similar role as separation oracles) were not
known. Moreover, in the context of polytopes, we obtain an adaptive
variant of our algorithm (see Algorithm~\ref{alg:dist_LMO_adaptive} and
Remark~\ref{rem:exact}) that uses additional linear programming
techniques and which is parameter-free. This variant
allows for \emph{exact} decision of whether there exists
$x \in P \cap Q$ (and computing it) or
deciding disjointness in number of required LMO calls
only a factor $4$ worse compared to the approximate version
as shown in Proposition~\ref{prop:adaptive}.
Moreover, if the polytopes are rational then
the logarithmic number of linear programs we need to solve
can be solved in polynomial time, so that the running time is
determined by the number of iterations of the
Algorithm~\ref{alg:dist_LMO_adaptive}.

An interesting corollary of our algorithms is that we can check whether
two polytopes (approximately or exactly) intersect \emph{without} an explicit
description of the polytopes by inequalities or vertices,
requiring only access via a linear minimization
oracle. Moreover, we obtain an algorithm that decides (non-)intersection for polytopes in \emph{finite time}, which is in stark contrast to the general case: Alternating projections converge in finitely many steps under
favorable circumstances, see \cite{behling2021infeasibility},
e.g., if \(P\) is a polytope and \(Q\) a hyperplane,
but there are polytopes \(P\) and \(Q\) for which they do not
converge in finitely many steps,
see e.g., \cite[Example~5.2]{behling2021infeasibility}.

\paragraph{Numerical Experiments}

We complement our theoretical results with numerical experiments in Section~\ref{sec:experiments}. While in the case where projections are easy, von Neumann's algorithm tends to be faster (as expected), in the case where projections are hard, ALM is significantly faster.

\paragraph{Follow-up Works}
Finally, we would like to mention that since our initial writing of this work several follow-up lines directly emerged from this work. One line of work considered \lmo-based splitting algorithms, with single-\lmo-call iterations for general optimization \citep{WP2023}, lazy block activation strategies to balance expensive vs cheap \lmo{}s \citep{BPW2024}, as well as strategies to combine the \lmo-based splitting with direct LP solving as done in \citet{HRBDP2025} as well as asynchronous activation strategies in \citet{SGDPW2024}. Alternating Linear Minimization (ALM) and its variants are now also standard algorithms in the \texttt{FrankWolfe.jl} package \citep{BDHKPTVW2024} and are part of the regular benchmark suite.

Our work also sparked interest in discrete geometry regarding the
minimal realizable distance of rational polytopes. Lower bounds were
known due to \citet{alon1997anti,graham1984anti}, and
\citet{DOPP2023,deza2025small,deza2025kissing3d,deza2026flat} showed
that these can be realized (up to small factors). This also shows that
our bounds presented here are effectively tight, both in iterations
vs.~distance, but also in terms of having the distance of the
polytopes appear in our complexities (see discussion in
\citet{DOPP2023}). Motivated by these works we also added a
complementary results in Remark~\ref{rem:adaptive_afw} that shows that
if one is willing to give up on alternations and incremental
activations one can achieve linear convergence in the position-related
measures with constants only depending on the polytopes separately and
not their interaction (i.e., whether they intersect or not). It
remains open whether in this case the pyramidal width constants of the
two polytopes have to be basically of the same worst-case order as in
the alternating case. 

We envision that our algorithm will have similar applications as
alternating projections, however additionally even in cases when
projections are hard or not available, as e.g., arising when
intersecting two compact convex sets with non-trivial
geometry. Well-known examples include the convex hull of
doubly-stochastic matrices (Birkhoff polytope), the matching polytope,
or the spectrahedron, where it is known that projections are much
slower than linear minimization (see e.g., \citet{CP2021}).

\subsubsection*{Outline}
In Section~\ref{sec:vonNeumann} we provide a quick recap of von
Neumann's alternating projection method, also called POCS,
with a short proof for
motivation. In Section~\ref{sec:alm} we then recall the Cyclic Block-Coordinate Conditional Gradient algorithm of \citet{beck2015cyclic} which is the basis for our alternating linear minimization algorithm, which we present in Section~\ref{sec:applications}.

\subsubsection*{Preliminaries}

In the following let \(\dist(X,Y)\) denote the Euclidean distance
between sets \(X\) and \(Y\),
and \(\dist(x, X)\) the distance between a point \(x\) and a set \(X\).
The \emph{Frank–Wolfe gap} of a differentiable function
$f \colon P \to \mathbb{R}$ at a point $x \in P$ is
\begin{equation}
  \label{eq:fwgap}
  \max_{v \in P} \innp{\nabla f(x)}{x - v}.
\end{equation}
Let \(\vertex(P)\) denote the set of vertices of a polytope \(P\).
Let \(\conv(X)\) denote the convex hull of a set \(X\) of points.

All other notation and notions are standard as to be found in \cite{CGFWSurvey2022}. 

\section{von Neumann's Alternating Projections}
\label{sec:vonNeumann}
We first briefly recall von Neumann's original algorithm (see
\cite{von1949rings})
and present a very elementary proof of its convergence.
While von Neumann's original proof applies only to linear subspaces,
we consider projections to convex sets;
in this setup the algorithm has also been called POCS
(Projections Onto Convex Sets) and
its convergence is well-known,
see e.g., \citet[Corollary~5.24]{bc17}.
Here we provide a likely well-known explicit convergence rate,
but we were not able to find a reference.
We provide a proof to help understanding later proofs.

We consider two closed convex sets $P, Q \subseteq \RR^n$ with
projectors $\Pi_P$ and $\Pi_Q$, respectively, where for a given $z \in
\RR^n$ we define $\Pi_P(z) \doteq \argmin_{x \in P} \norm{x - z}^2$
and similarly for $\Pi_Q(z)$, i.e., $\Pi_P(z)$ and $\Pi_Q(z)$ are the
projections of $z$ onto $P$ and $Q$, respectively, under the $2$-norm.
We shall use the following inequality for projection
\cite[proof of Theorem~E.9.3.0.1]{convexBook}:
\begin{equation}
  \label{eq:projection}
  \norm{x - y}^{2}
  \geq
  \norm{\Pi_{P}(x) - \Pi_{P}(y)}^{2}
  +
  \norm{x - \Pi_{P}(x) - y + \Pi_{P}(y)}^{2}
\end{equation}
generalizing  for \(x \in P\) the well-known inequality \(\norm{x - y}^{2} \geq
\norm{x - \Pi_{P}(y)}^{2} + \norm{\Pi_{P}(y) - y}^{2}\),
which is equivalent to \(\innp{x - \Pi_{P}(y)}{\Pi_{P}(y) - y} \geq 0\).
If \(P\) and \(Q\) are compact then by continuity of distance
there is \(x^{*} \in P\) and \(y^{*} \in Q\) with minimal distance,
obviously \(y^{*} = \Pi_{Q}(x^{*})\) and \(x^{*} = \Pi_{P}(y^{*})\).
We call \(x^{*} -y^{*}\) the \emph{distance vector}
between \(P\) and \(Q\),
which is easily seen by the above inequality to be independent of
the choice of \(x^{*}\) and \(y^{*}\).
Indeed, for any \(x \in P\) and \(y \in Q\),
by the above
we have \(\innp{x - x^{*}}{x^{*} - y^{*}} \geq 0\)
\begin{equation}
  \label{eq:min-dist}
 \begin{split}
  \innp{x - y}{x^{*} - y^{*}}
  &
  =
  \innp{x - x^{*}}{x^{*} - y^{*}}
  + \innp{x^{*} - y^{*}}{x^{*} - y^{*}}
  + \innp{y^{*} - y}{x^{*} - y^{*}}
  \\
  &
  \geq
  0 + \dist(P, Q)^{2} + 0 = \dist(P, Q)^{2}
  .
 \end{split}
\end{equation}
Thus if \(x \in P\) and \(y \in Q\) also have minimal distance,
i.e., \(\norm{x-y} = \dist(P, Q) = \norm{x^{*} - y^{*}}\),
then necessarily \(x - y = x^{*} - y^{*}\).

\begin{algorithm}[H]
  \caption{POCS, von Neumann's Alternating Projections}
	\label{alg:vnap}  
	\begin{algorithmic}[1]
		\Require Point $y_{0} \in \RR^n$, $\Pi_P$ projector onto $P \subseteq \RR^n$ and $\Pi_Q$ projector onto $Q \subseteq \RR^n$. 
    \Ensure Iterates $x_1,y_1 \dotsc \in \RR^n$
		\hrulealg
                \For{$t=0$ \To \dots}
		\State $x_{t+1} \leftarrow \Pi_P(y_{t})$
		\State $y_{t+1} \leftarrow \Pi_Q(x_{t+1})$
		\EndFor
	\end{algorithmic}
\end{algorithm}

\begin{figure}
  \centering
  \begin{tikzpicture}[scale=.3, >=stealth]
    \draw[thick] (0, 5) -- (4, 9) -- node[above]{\(Q\)} (-5, 5)
    --cycle;
    \draw[thick] (-4, 0) -- (6, 0) -- (12, -2) -- node[below]{\(P\)}
    cycle;
    \fill (10, 2) circle[radius=7pt];
    \draw[postaction={on each segment={mid arrow}}]
    (10, 2) node[right]{\(y_{0}\)}
    -- (9,-1) node[right, yshift=.7ex]{\(x_{1}\)}
    -- (2, 7) node[right=1ex]{\(y_{1}\)}
    -- (2, 0) node[above right]{\(x_{2}\)}
    -- (0, 5) node[above, xshift=-1ex]{\(y_{2}\)}
    -- (0, 0) node[above left]{\(x_{3}\)};
    \draw[->, dashed] (-3,5) -- node[left]{\(d\)} (-3,0);
  \end{tikzpicture}
  \caption{\label{fig:alt-proj}
    Alternating projections between two triangles \(P\)
    and \(Q\) with distance vector \(d\).
    In this simple example the algorithm alternated between
    a pair of minimal distance points (\(y_{2}\) and \(x_{3}\)) after
    finitely many projections.}
\end{figure}

\begin{proposition}
	\label{prop:vnAlt}
  Let \(P\) and \(Q\) be compact convex sets and
  let $x_1,y_1 \dotsc, x_T, y_T \in \RR^n$ be the sequence of iterates
  of Algorithm~\ref{alg:vnap}.
  Let \(Q_{\min}\) be the sets of points of \(Q\) with minimal
  distance to \(P\), and \(d\) the distance vector
  between \(Q\) and \(P\),
  i.e., \(d = z - \Pi_{P}(z)\) for all \(z \in Q_{\min}\).
  (See Figure~\ref{fig:alt-proj} for an example.)
  Then the iterates converge: \(x_{t} \to x\) and \(y_{t} \to y\)
  to some \(x \in P\) and \(y \in Q\) with \(x - y = d\).
  Moreover, for \(T \geq 1\)
  \begin{equation}
    \label{eq:alt-proj-conv}
    \min_{t=0, \dots, T-1} \norm{y_{t} - x_{t+1} -d}^{2}
    + \norm{x_{t+1} - y_{t+1} + d}^{2}
    \leq
    \frac{1}{T} \sum_{t = 0}^{T-1} \left(\norm{y_{t} - x_{t+1} - d}^{2}
      + \norm{x_{t+1} - y_{t+1} + d}^{2} \right)
    \leq
    \frac{\dist(y_0, Q_{\min})^{2}}{T}
    .
  \end{equation}
  In particular, if \(P\) and \(Q\) intersect
  (i.e., \(Q_{\min} = P \cap Q\) and \(d=0\))
  then
  \begin{equation}
    \label{eq:alt-proj-conv-intersect}
    \norm{x_{T} - y_{T}}^{2}
    \leq
    \frac{1}{T} \sum_{t = 0}^{T-1} \left(\norm{y_{t} - x_{t+1}}^{2}
      + \norm{x_{t+1} - y_{t+1}}^{2} \right)
    \leq
    \frac{\dist(y_0, P \cap Q)^{2}}{T}
    .
  \end{equation}
\begin{proof}
Let \(z_{2} \in Q_{\min}\) and \(z_{1} = \Pi_{P}(z_{2}) \in P\).
Obviously, \(\norm{z_{1} - z_{2}} = \dist(P, Q)\),
and hence \(z_{2} = \Pi_{Q}(z_{1})\) and \(d = z_{2} - z_{1}\).

We estimate the summands in the claim
via the projection inequality~\eqref{eq:projection}:
\begin{subequations}
  \label{eq:Neuman-decrease}
  \begin{align}
    \norm{x_{t+1} - y_{t+1} + d}^{2} =
    \norm{x_{t+1} - y_{t+1} - z_{1} + z_{2}}^{2} & \leq
    \norm{x_{t+1} - z_{1}}^{2} - \norm{y_{t+1} - z_{2}}^{2}
    \\
    \norm{y_{t} - x_{t+1} - d}^{2} =
    \norm{y_{t} - x_{t+1} - z_{2} + z_{1}}^{2} & \leq
    \norm{y_{t} - z_{2}}^{2} - \norm{x_{t+1} - z_{1}}^{2} .
  \end{align}
\end{subequations}
Summing up provides a telescope sum:
\begin{multline*}
  \sum_{t=0}^{T-1} \left(
    \norm{x_{t+1} - y_{t+1} + d}^{2}
    + \norm{y_{t} - x_{t+1} - d}^{2}
  \right)
  \\
  \leq
  \sum_{t=0}^{T-1} \left[
    (\norm{x_{t+1} - z_{1}}^{2} - \norm{y_{t+1} - z_{2}}^{2})
    +
    (\norm{y_{t} - z_{2}}^{2} - \norm{x_{t+1} - z_{1}}^{2})
  \right]
  \\
  =
  \norm{y_{0} - z_{2}}^{2} - \norm{y_{T} - z_{2}}^{2}
  \leq
  \norm{y_{0} - z_{2}}^{2}
  .
\end{multline*}
Minimizing over \(z_{2}\) proves the second inequality of
Equation~\eqref{eq:alt-proj-conv}.
The first inequality of Equation~\eqref{eq:alt-proj-conv}
is obvious.

From Equation~\eqref{eq:alt-proj-conv} it follows that
\(x_{t} - y_{t}\) converges to \(-d\).
Next we prove that the individual sequences \(x_{t}\) and \(y_{t}\)
converge.
As \(Q\) is compact, the sequence \(y_{t}\)
has some accumulation point \(y\).
Then \(x \doteq y-d\) is an accumulation point of \(x_{t}\) and
thus \(x \in P\) and \(y \in Q\) with \(x - y = -d\),
hence \(y \in Q_{\min}\).

By Equations~\eqref{eq:Neuman-decrease}, we have
\(\norm{y_{t} - z_{2}}^{2} \geq \norm{x_{t+1} - \Pi_{P}(z_{2})}^{2} \geq
\norm{y_{t+1} - z_{2}}^{2}\)
for every \(z_{2} \in Q_{\min}\), in particular,
\(\norm{y_{t} - z_{2}}\) is a decreasing sequence.
This means for the accumulation point \(z_{2} \doteq y\)
that \(\norm{y_{t} - y}\) must converge to \(0\),
i.e., \(y_{t}\) converges to \(y\).
Therefore \(x_{t}\) converges to \(y-d=x\).

Finally, when \(P\) and \(Q\) intersect, i.e., \(d = 0\)
then as projections select minimal distance points,
the distances are decreasing:
\(\norm{x_{t} - y_{t}}^{2} \geq \norm{x_{t+1} - y_{t}}^{2} \geq
\norm{x_{t+1} - y_{t+1}}^{2}\),
and therefore the minimum in \eqref{eq:alt-proj-conv} simplifies.
\end{proof}
\end{proposition}

\section{The Cyclic Block-Coordinate Conditional Gradient algorithm}
\label{sec:alm}

Here we recall the Cyclic
Block-Coordinate Conditional Gradient algorithm (CBCG) of
\citet{beck2015cyclic} as Algorithm~\ref{alg:alt_lin_min}, 
using the notation of its predecessor from
\citet{lacoste2013block}.
The algorithm solves the optimization problem:
\begin{equation}
	\label{eq:optSolve}
	\min_{(x_{0},\dotsc,x_{k-1}) \in P_{0} \times \dotsm \times P_{k-1}}
	f(x_{0}, \dotsc, x_{k-1}),
\end{equation}
where $f$ is a convex function and the \(P_{i}\) are convex sets.
This problem can be also solved with the original Frank–Wolfe algorithm
\citep{frank1956algorithm}, however, for this paper
the cyclic variant is more suitable.
CBCG makes a single linear minimization in only one block \(P_{i}\)
in every iteration,
and blocks are selected in a simple cyclic order.

Let \(s_{[i]}\) denote the vector equal to \(s\) in \(P_{i}\)
and to \(0\) elsewhere, i.e.,
it is the natural generalization of the coordinate vectors. We now state the convergence rate of the Cyclic Block-Coordinate
Conditional Gradient algorithm from \citet{beck2015cyclic} with a
slight improvement in the constants for the agnostic step rule:
the denominator is improved from \(t + 2\) to \(t + 1\), and the
constant factor \(6.75\) for the dual gap is a bit smaller.

\begin{algorithm}
  \caption{Cyclic Block-Coordinate Conditional Gradient algorithm \citep{beck2015cyclic}}
	\label{alg:alt_lin_min}  
	\begin{algorithmic}[1]
          \Require Points $x^{0}_{i} \in P_{i}$,
          LMO for $P_{i} \subseteq \RR^{n_{i}}$, \(i=0, \dotsc, k-1\)
          and $0 < \gamma_0, \dots, \gamma_t, \dotsc \leq 1$.
                \Ensure Iterates $x^{1}, \dotsc \in P_{0} \times
                \dotsm \times P_{k-1}$
		\hrulealg
                \For{$t=0$ \To $\dotsc$}
                \State \(i \gets t \bmod k\)
                \State\label{alg:fw_vertex^P}
                \(v^{t} \gets \argmin_{x \in P_{i}}
                \innp{\nabla_{P_{i}} f(x^{t})}{x}\)
                \State\label{alg:progress}\(x^{t+1} \gets x^{t}
                + \gamma_{t} (v^{t} - x^{t}_{i})_{[i]}\)
		\EndFor
	\end{algorithmic}
\end{algorithm}

\begin{theorem}[{Convergence of Cyclic Block-Coordinate
    Conditional Gradient algorithm \citep[cf
    Theorem~4.5]{beck2015cyclic}}]
  \label{thm:almo}
  Let $P_{i}$ be a compact convex set with diameter $D_{i}$
  for \(i=0, 1, \dotsc, k-1\).
  Let $f \colon P_{0} \times \dotsm \times P_{k-1} \to \RR$ be
  an \(L\)-smooth convex function,
  which is moreover partially \(L_{i}\)-smooth in \(P_{i}\)
  for \(i=0, 1, \dotsc,
  k-1\).
  Let \(D \doteq \sqrt{\sum_{i=0}^{k-1} D_{i}^{2}}\) be the diameter
  of the feasible region \(P_{0} \times \dotsm \times P_{k-1}\).
  Then the iterates of Algorithm~\ref{alg:alt_lin_min} with the choice
  \(\gamma_{t} = \frac{2}{\lfloor t / k \rfloor + 2}\) satisfy
  \begin{align*}
    f(x^{k t}) - f(x^{*}) &\leq \frac{2}{t + 2}
    \left(
      \sum_{i=0}^{k-1} \frac{L_{i} D_{i}^{2}}{2}
      + 2 L D \sum_{i=0}^{k-1} D_{i}
    \right),
    \\
    \min_{1 \leq t \leq T}
    \max_{y \in P_{0} \times \dots \times P_{k-1}}
    \innp{\nabla f(x^{k t})}{x^{k t} - y}
    &
    \leq
    \frac{6.75}{T + 2}
    \left(
      \sum_{i=0}^{k-1} \frac{L_{i} D_{i}^{2}}{2}
      + 2 L D \sum_{i=0}^{k-1} D_{i}
    \right)
    .
  \end{align*}
  where $x^*$ is a minimizer of \(f\).
\begin{proof}
By \citet[Lemma~4.3]{beck2015cyclic}
\begin{equation}
  \label{eq:alt-progress}
  f(x^{k t + k}) - f(x^{k t})
  - \gamma_{k t} \innp{\nabla f(x^{k t})}{u^{k t} - x^{k t}}
  \leq
  \gamma_{k t}^{2}
  \left(
    \sum_{i=0}^{k-1} \frac{L_{i} D_{i}^{2}}{2}
    + 2 L D \sum_{i=0}^{k-1} D_{i}
 \right)
  .
\end{equation}
This differs from the standard recursion for the vanilla Frank–Wolfe
algorithm only by replacing the smoothness constant with the
parenthesized expression on the right multiplied by \(2\).
Therefore the rest of the proof is standard, see e.g.,
\cite{jaggi2013revisiting} or \cite{CGFWSurvey2022}.
\end{proof}
\end{theorem}

\begin{theorem}[{Convergence of Cyclic Block-Coordinate
    Conditional Gradient algorithm \citep[Theorem~4.13]{beck2015cyclic}}]
  \label{thm:almo-short}
  Let $P_{i}$ be a compact convex set with diameter $D_{i}$
  for \(i=0, 1, \dotsc, k-1\).
  Let $f \colon P_{0} \times \dotsm \times P_{k-1} \to \RR$ be
  an \(L\)-smooth convex function,
  which is moreover partially \(L_{i}\)-smooth
  and partially \(G_{i}\)-Lipschitz in \(P_{i}\)
  for \(i=0, 1, \dotsc,
  k-1\).
  Let \(D \doteq \sqrt{\sum_{i=0}^{k-1} D_{i}^{2}}\) be the diameter
  of the feasible region \(P_{0} \times \dotsm \times P_{k-1}\).
  Then the iterates of Algorithm~\ref{alg:alt_lin_min} with the
  short step-size rule
  \(\gamma_{t} = \min \{\innp{\nabla_{P_{i}} f(x^{t})}{x^{t}_{i} - v^{t}}
  / (L_{i} \norm{x^{t}_{i} - v^{t}}^{2}), 1\}\) satisfy
  \begin{align*}
    f(x^{k t}) - f(x^{*}) &\leq
    \frac{4 k}{t + 4}
    \left(
      \max_{i = 0, \dots, k-1} \{L_{i} D_{i}^{2}, G_{i} D_{i}\}
      +
      \frac{k L^{2} D^{2}}{\min_{i = 0, \dots, k-1} L_{i}}
    \right),
    \\
    \min_{1 \leq t \leq T}
    \max_{y \in P_{0} \times \dots \times P_{k-1}}
    \innp{\nabla f(x^{k t})}{x^{k t} - y}
    &
    \leq
    \frac{8 k}{T + 4}
    \left(
      \max_{i = 0, \dots, k-1} \{L_{i} D_{i}^{2}, G_{i} D_{i}\}
      +
      \frac{k L^{2} D^{2}}{\min_{i = 0, \dots, k-1} L_{i}}
    \right)
    .
  \end{align*}
  where $x^*$ is a minimizer of \(f\).
\end{theorem}

Note that typically the short step-size rule used in
Theorem~\ref{thm:almo-short} performs much better in actual
computations than the agnostic step-size rule used in
Theorem~\ref{thm:almo}. Therefore we will state our results for both
choices in the remainder.

\section{Alternating Linear Minimizations}
\label{sec:applications}

We are returning to our original problem of interest:
Let $P$ and $Q$ be compact convex sets contained in the same ambient space
$\RR^n$, and the goal is to find points in them with minimal distance.
In contrast to von Neumann's approach, we will assume that we have
only access to linear minimization oracles (LMOs) for the compact
convex sets $P \subseteq \RR^n$ and $Q \subseteq \RR^n$, i.e., given a
vector $c \in \RR^n$ the oracles return:
\[
v \leftarrow \argmin_{x \in P} \innp{c}{x} \qquad \text{and} \qquad w
\leftarrow \argmin_{y \in Q} \innp{c}{y}.
\]

We adapt Algorithm~\ref{alg:alt_lin_min} to find points in \(P\) and \(Q\)
with approximate minimal distance, using the objective function
\(f(x,y) = \norm{x-y}^{2}\) over \(P \times Q\).

\begin{algorithm}[H]
  \caption{Alternating Linear Minimizations (ALM)}
  \label{alg:dist_LMO}
  \begin{algorithmic}[1]
    \Require Points \(x_{0} \in P\), $y_{0} \in Q$,
      LMO over \(P, Q \subseteq \RR^{n}\)
    \Ensure Iterates $x_1,y_1 \dotsc \in \RR^n$
    \hrulealg
    \For{$t=0$ \To \dots}
      \State \(u_{t} \gets \argmin_{x \in P} \innp{x_{t} - y_{t}}{x}\)
      \State \(x_{t+1} \gets
        x_{t} + \gamma_{t,1} \cdot (u_{t} - x_{t})\)
      \State \(v_{t} \gets \argmin_{y \in Q} \innp{y_{t} - x_{t+1}}{y}\)
      \State \(y_{t+1} \gets
        y_{t} + \gamma_{t,2} \cdot (v_{t} - y_{t})\)
    \EndFor
  \end{algorithmic}
\end{algorithm}

Theorems~\ref{thm:almo} and~\ref{thm:almo-short}
provide the following guarantees:

\begin{proposition}[Intersection of two sets]
	\label{prop:setIntersection}
  Let $P$ and $Q$ be compact convex sets.
  Then Algorithm~\ref{alg:dist_LMO}
  with \(\gamma_{t,1} = \gamma_{t,2} \doteq \frac{2}{t+2}\) (agnostic step-size rule)
  generates iterates $z_t \doteq \frac{1}{2} (x_t + y_t)$, such that
  \begin{align}
    \max\{\dist(z_t, P)^2, \dist(z_t, Q)^2\}
    \leq
    \frac{\norm{x_{t} - y_{t}}^{2}}{4}
    &
    \leq \frac{(1 + 2 \sqrt{2}) (D_{P}^{2} + D_{Q}^{2})}{t+2}
    + \frac{\dist(P,Q)^{2}}{4}
    \\
    \label{eq:P-cap-Q-dual}
    \min_{1 \leq t \leq T}
    \max_{x \in P, y \in Q}
    \norm{x_{t} - y_{t}}^{2} - \innp{x_{t} - y_{t}}{x - y}
    &
    \leq
    \frac{6.75 (1 + 2 \sqrt{2})}{T + 2}
    (D_{P}^{2} + D_{Q}^{2})
    .
  \end{align}

  With the short step-size rule: \(\gamma_{t,1} \doteq \min\{
  \innp{x_{t} - y_{t}}{x_{t} - u_{t}} / \norm{x_{t} - u_{t}}^{2}, 1\}\)
  and similarly
  \(\gamma_{t,2} \doteq \min\{
  \innp{y_{t} - x_{t+1}}{y_{t} - v_{t}} / \norm{y_{t} - v_{t}}^{2},
  1\}\)
  the iterates of Algorithm~\ref{alg:dist_LMO} satisfy
  \begin{align}
    \max\{\dist(z_t, P)^2, \dist(z_t, Q)^2\}
    \leq
    \frac{\norm{x_{t} - y_{t}}^{2}}{4}
    &
    \leq
    \frac{4 c}{t+4}
    + \frac{\dist(P,Q)^{2}}{4}
    \\
    \min_{1 \leq t \leq T}
    \max_{x \in P, y \in Q}
    \norm{x_{t} - y_{t}}^{2} - \innp{x_{t} - y_{t}}{x - y}
    &
    \leq
    \frac{16 c}{T+4}
  \end{align}
  where
  \begin{equation}
    \label{eq:1}
    c \doteq
    (D_{P} + D_{Q} + \dist(P,Q)) \cdot \max \{D_{P}, D_{Q}\}
    + 2 (D_{P}^{2} + D_{Q}^{2})
    .
  \end{equation}
\begin{proof}
First note that $\dist(z_t, P)^2 \leq \norm{x_t-z_t}^2 = \norm{x_t -
  \frac{1}{2} (x_t + y_t)}^2 = \norm{x_t - y_t}^2 / 4$ and similarly
$\dist(z_t, Q) \leq \norm{x_t - y_t} / 2$.
Combining these with Theorems~\ref{thm:almo} and~\ref{thm:almo-short},
we obtain
the claim.
\end{proof}
\end{proposition}

It is informative to compare the rate obtained in Proposition~\ref{prop:setIntersection} with the rate in Proposition~\ref{prop:vnAlt}. 

\begin{remark}[Comparison to von Neumann's alternating projection algorithm]
  For simplicity let us consider the case where
  $P \cap Q \neq \emptyset$.
  Via Proposition~\ref{prop:vnAlt} we then obtain the following
  guarantee for von Neumann's alternating projection algorithm:
  \begin{equation}
    \norm{x_{T} - y_{T}}^{2}
	\leq
	\frac{\dist(y_0, P \cap Q)^{2}}{T}
	.
\end{equation}
With $z_{t} \doteq \frac{1}{2} (x_{t} + y_{t})$, we have
\begin{equation*}
  \max\{\dist(z_{t+1}, P)^2, \dist(z_{t+1}, Q)^2\}
  \leq \norm{x_{t+1} - z_{t+1}}^2 + \norm{y_{t+1} - z_{t+1}}^2
  = \frac{1}{2} \norm{x_{t+1} - y_{t+1}}^2.
\end{equation*}
Combining with the above and estimating generously we obtain for von Neumann's alternating projection method:
\begin{equation}
	\label{eq:simplifiedVNg}
	\min_{t=0, \dots, T-1} \max\{\dist(z_{t+1}, P)^2, \dist(z_{t+1}, Q)^2\} \leq \frac{\dist(y_0, P \cap Q)^{2}}{T}.
\end{equation}
In contrast to this, the LMO-based approach guarantees via Proposition~\ref{prop:setIntersection} after $T$ many iterations:
\[
\max\{\dist(z_T, P)^2, \dist(z_T, Q)^2\}
\leq \frac{(1 + 2 \sqrt{2}) (D_{P}^{2} + D_{Q}^{2})}{T+2}
.
\]	
So both approaches essentially guarantee an $O(1/T)$
rate of convergence, however the constants are quite different and not
necessarily comparable, depending on the regime. The vanilla
Frank–Wolfe algorithm applied directly to the function
$f(x,y) = \norm{x-y}^2$
would essentially lead to the same convergence guarantee as
Algorithm~\ref{alg:dist_LMO} (up to a constant factor $\sqrt{2}$),
however the processing would not be alternating between the blocks
anymore.

Finally, running von Neumann's original algorithm with a näive
simulation of the projection operator via the vanilla Frank–Wolfe
method (which only requires LMOs) would roughly require an additional
number of $\Omega(t)$
LMO calls per iteration. This assumes the target accuracy is at least
$O(1/t)$
which might be required to solve the projection problem with
reasonably high accuracy.  Thus the overall convergence rate in the
number of LMO calls would be roughly $O(1/\sqrt{T})$,
i.e., much slower than our Algorithm~\ref{alg:dist_LMO}.
\end{remark}

\subsection{Emptiness of intersection}
\label{sec:empty-intersection}

We can also use the same approach to certify that the intersection of
two sets $P$ and $Q$ is empty.
To this end, we simply run the algorithm until the
function value $f(x_t,y_t)$ is strictly larger than the primal or dual gap.

\begin{corollary}[Certifying that $P \cap Q = \emptyset$]
  \label{cor:disjoint}
  Let $P, Q$ be compact convex sets with diameters \(D_{p}\) and
  \(D_{Q}\), respectively.
  If some iterates of Algorithm~\ref{alg:dist_LMO}
  using \(\gamma_{t,1} = \gamma_{t,2} = 2/(t+2)\)
  satisfy
  \begin{equation*}
    \norm{x_{t} - y_{t}}^{2}
    > \frac{4 (1 + 2 \sqrt{2}) (D_{P}^{2} + D_{Q}^{2})}{t+2}
  \end{equation*}
  then \(P\) and \(Q\) are disjoint.
  If \(P\) and \(Q\) are disjoint then the above condition
  is satisfied after at most
  \[
    \frac{8 (1 + 2 \sqrt{2}) (D_{P}^{2} + D_{Q}^{2})}{\dist(P,Q)^{2}}
  \]
  block-LMO calls.

  Under the short step rule, \(P\) and \(Q\) are disjoint if for some
  iterate \(t\)
  \begin{equation}
    \label{eq:short-disjoint}
    \norm{x_{t} - y_{t}}^{2}
    >
    \frac{16}{t+4}
    \left(
      (D_{P} + D_{Q}) \cdot \max \{D_{P}, D_{Q}\}
      + 2 (D_{P}^{2} + D_{Q}^{2})
    \right)
    .
  \end{equation}
  If \(P\) and \(Q\) are disjoint then this condition
  is satisfied after at most the following number of block-LMO calls
  \begin{equation}
    \frac{32}{(t+4) \dist(P,Q)^{2}}
    \left(
      (D_{P} + D_{Q}) \cdot \max \{D_{P}, D_{Q}\}
      + 2 (D_{P}^{2} + D_{Q}^{2})
    \right)
    .
  \end{equation}
\begin{proof}
We prove only the claims for the agnostic step step size rule
\(\gamma_{t,1} = \gamma_{t,2} = 2/(t+2)\)
as the proof for the short step rule is similar.
We use
	the following chain of inequalities,
	the first inequality
	by Proposition~\ref{prop:setIntersection},
	the second inequality because \(x_{t} \in P\) and \(y_{t} \in Q\):
	\[
	\dist(P,Q)^{2}
	\geq
  \norm{x_{t} - y_{t}}^{2}
  - \frac{4 (1 + 2 \sqrt{2}) (D_{P}^{2} + D_{Q}^{2})}{t+2}
	\geq
  \dist(P,Q)^{2} - \frac{4 (1 + 2 \sqrt{2}) (D_{P}^{2} + D_{Q}^{2})}{t+2}
	.
	\]
	By the first inequality,
  \(\norm{x_{t} - y_{t}}^{2}
  > \frac{4 (1 + 2 \sqrt{2}) (D_{P}^{2} + D_{Q}^{2})}{t+2}\)
	obviously implies \(\dist(P, Q) > 0\), i.e., that \(P\) and \(Q\) are
	disjoint.
	By the second inequality,
  \(\norm{x_{t} - y_{t}}^{2}
  > \frac{4 (1 + 2 \sqrt{2}) (D_{P}^{2} + D_{Q}^{2})}{t+2}\)
	for \(t >  4 (D_{p} + D_{Q})^{2} / \dist(P,Q)^{2}\).
	As every iteration of Algorithm~\ref{alg:dist_LMO}
	costs two block-LMO calls, this provides the upper bound
  \(8 (1 + 2 \sqrt{2}) (D_{P}^{2} + D_{Q}^{2}) / \dist(P,Q)^{2}\)
  on the number of block-LMO calls.
\end{proof}
\end{corollary}

With a slightly worse constant factor we can also certify and test for disjointness without knowing $D_P$ and $D_Q$ explicitly. 

\begin{corollary}[Certifying $P \cap Q = \emptyset$
  without diameters]
  \label{cor:disjoint-par-free}
  Let $P, Q$ be disjoint compact convex sets
  with diameters \(D_{p}\) and \(D_{Q}\), respectively.
  Then executing Algorithm~\ref{alg:dist_LMO}, after at most
  \begin{gather*}
    \frac{13.5 (1 + 2 \sqrt{2}) (D_{P}^{2} + D_{Q}^{2})}{\dist(P,Q)^{2}}
    \qquad
    \text{with step sizes \(\gamma_{t,1} = \gamma_{t,2} = 2/(t+2)\)}
    \\
    \frac{32}{\dist(P,Q)^{2}}
    \left(
      (D_{P} + D_{Q} + \dist(P,Q)) \cdot \max \{D_{P}, D_{Q}\}
      + 2 (D_{P}^{2} + D_{Q}^{2})
    \right)
    \qquad
    \text{with short step rule}
  \end{gather*}
  block-LMO calls, some (of the already seen!) iteration \(t\) provides
  the following certificate for disjointness, which does not require explicit bounds on $D_P$ and $D_Q$ (in contrast to the above):
  \begin{equation}
    \label{eq:P-Q-sep}
    \min_{x \in P, y \in Q} \innp{x_{t} - y_{t}}{x - y} > 0
    .
  \end{equation}
  Moreover, this inequality is guaranteed to hold for every iteration
  \(t > 4 (1 + 2 \sqrt{2}) (D_{P}^{2} + D_{Q}^{2}) (D_{P} + D_{Q})^{2}
  / \dist(P, Q)^{4}\),
  i.e., after
  \(8 (1 + 2 \sqrt{2}) (D_{P}^{2} + D_{Q}^{2}) (D_{P} + D_{Q})^{2}
  / \dist(P, Q)^{4}\)
  block-LMO calls
  under the agnostic step size rule
  \(\gamma_{t,1} = \gamma_{t,2} = 2/(t+2)\),
  and after
  \(32
    [
      (D_{P} + D_{Q} + \dist(P, Q)) \cdot \max \{D_{P}, D_{Q}\}
      + 2 (D_{P}^{2} + D_{Q}^{2})
    ] (D_{P} + D_{Q})^{2} / \dist(P,Q)^{4}\)
  block-LMO calls
  under the short step rule.
\begin{proof}
Clearly, \eqref{eq:P-Q-sep} implies disjointness of \(P\)
and \(Q\).
Let us assume that
\(\min_{x \in P, y \in Q} \innp{x_{t} - y_{t}}{x - y} \leq 0\)
for all \(1 \leq t \leq T\).
By \eqref{eq:P-cap-Q-dual} of Proposition~\ref{prop:setIntersection},
  \begin{multline*}
    \dist(P, Q)^{2}
    \leq
    \min_{1 \leq t \leq T}
    \norm{x_{t} - y_{t}}^{2}
    \leq
    \min_{1 \leq t \leq T}
    \max_{x \in P, y \in Q}
    \norm{x_{t} - y_{t}}^{2} - \innp{x_{t} - y_{t}}{x - y}
    \leq
    \frac{6.75 (1 + 2 \sqrt{2})}{T + 2}
    (D_{P}^{2} + D_{Q}^{2})
    .
  \end{multline*}
  Therefore
  \(T + 2 \leq \frac{6.75 (1 + 2 \sqrt{2})}{\dist(P, Q)^{2}}
  (D_{P}^{2} + D_{Q}^{2})\).
  Thus, Inequality~\eqref{eq:P-Q-sep}
  is satisfied by some of the first
  \(\frac{6.75 (1 + 2 \sqrt{2})}{\dist(P, Q)^{2}}
  (D_{P}^{2} + D_{Q}^{2})\) iterations.
  The upper bound on the block LMO calls to satisfy the inequality
  is twice this number since every iteration costs two block LMO calls
  as above.
\end{proof}
\end{corollary}

\begin{remark}[$V$-representation vs.~$H$-representation vs.~implicit description]
  Given two polytopes $P$ and $Q$ in $H$-representation,
  i.e., each given as a set of defining inequalities, then it is
well-known to decide whether $P \cap Q = \emptyset$
  and compute $x \in P \cap Q$
  via Farkas' lemma and linear programming techniques.

  If the polytopes are given in $V$-representation
  (i.e., as finite set whose convex hull is the polytope), then we can
  resolve the question similarly: Let $P = \conv(U)$ and $Q =
  \conv(W)$.
  We consider the extended formulations arising from $U$ and
  $W$ and check the feasibility of the linear program
  \[
  \left\{(\lambda, \kappa) : \sum_{u \in U} \lambda_u u = \sum_{w \in W} \kappa_w w, \sum_{u \in U} \lambda_u = \sum_{w \in W} \kappa_w = 1, \lambda, \kappa \geq 0\right\}.
  \]
  The mixed case can be decided similarly. 
  
  If $P$ and $Q$ are only given implicitly via linear optimization oracles, it is \emph{a priori} not obvious how to decide
  whether the polytopes intersect. However, we can use 
  Algorithm~\ref{alg:dist_LMO} to, at least,
  approximately decide this question and use Remark~\ref{rem:exact} to make the test exact.
\end{remark}

\subsection{Adaptive variant for polytopes}
\label{sec:adaptive-alm-polytope}

The careful reader will have realized that Algorithm~\ref{alg:dist_LMO} only provides an approximate answer as to whether $x \in P \cap Q$ exists and in fact, it only decides whether there exists a point $x \in \RR^n$ that is close to both $P$ and $Q$. For the case where $P$ and $Q$ are polytopes, we will now show that there exists $\varepsilon_{PQ} > 0$ depending on $P$ and $Q$, such that if we obtain iterates $x_t$ and $y_t$  via Algorithm~\ref{alg:dist_LMO} with $\norm{x_t - y_t} \leq \varepsilon_{PQ}$, then we can compute an actual point $x \in P \cap Q$. 

To this end we start with a simple observation. 

\begin{observation}
	\label{obs:intersect}
Let $P, Q \subseteq \RR^n$ be polytopes. There exists
$\varepsilon_{PQ} > 0$, so that for all $U \subseteq \vertex(P),
V \subseteq \vertex(Q)$ with
$\dist(\conv(U), \conv(V)) < \varepsilon_{PQ}$,
it holds $\conv(U) \cap \conv(V) \neq \emptyset$.
\begin{proof}
This easily follows from
polytopes having only a finite number of vertices,
so that the simple definition provides a positive number:
\[\varepsilon_{PQ} \coloneqq \min\{\dist(\conv(U), \conv(V)) : U
  \subseteq \vertex(P), V \subseteq \vertex(Q), \conv(U) \cap \conv(V)
  = \emptyset\}.
\qedhere\]
\end{proof}
\end{observation}

With this observation we obtain now:

\begin{remark}[Recovery of $x \in P \cap Q$ by linear programming]
	\label{rem:exact}
  Once we obtain $x_t$ and $y_t$
  with $\norm{x_t - y_t} < \varepsilon_{PQ}$
  via Algorithm~\ref{alg:dist_LMO}, we can recover $x \in P \cap Q$
  via linear programming. To this end let $U \subseteq \vertex(P)$
  be all extreme points returned by the LMO for $P$
  throughout the execution of Algorithm~\ref{alg:dist_LMO} and define
  $V \subseteq \vertex(Q)$
  accordingly. From Observation~\ref{obs:intersect} it follows that
  $\conv(U) \cap \conv(V) \neq \emptyset$
  and we can solve the linear feasibility program:
\begin{subequations}
  \label{eq:LMO-point}
  \begin{align}
    \sum_{u \in U} \lambda_u u & = \sum_{v \in V} \kappa_u v \\
    \sum_{u \in U} \lambda_u & = 1 \\
    \sum_{v \in V} \kappa_u & = 1 \\
    \lambda_u & \geq 0 &\forall u \in U \\
    \kappa_v & \geq 0 &\forall v \in V.
  \end{align}
\end{subequations}
This linear feasibility problem is guaranteed to be feasible and the obtained solution satisfies $x \coloneqq \sum_{u \in U} \lambda_u u = \sum_{v \in V} \kappa_u v \in P \cap Q$.

If $\varepsilon_{PQ}$ is not known ahead of time, which will usually be the case, we can simply try to solve the above linear program whenever the distance $\norm{x_t - y_t}$ is halved. This leads to a logarithmic overhead of about $O(\log 1/\varepsilon_{PQ})$ solved linear programs.
\end{remark}

In fact we can turn Algorithm~\ref{alg:dist_LMO} into an exact
adaptive algorithm (Algorithm~\ref{alg:dist_LMO_adaptive}) for
polytopes not requiring knowledge of any parameters such as
$\varepsilon_{PQ}$ and diameters at the cost of an additive
logarithmic overhead of LMO calls.
The algorithm can be also run for more general regions, however then
the guarantees are not as clean.

\begin{algorithm}
	\caption{Alternating Linear Minimizations (ALM) (adaptive variant)}
	\label{alg:dist_LMO_adaptive}
	\begin{algorithmic}[1]
		\Require Points \(x_{0} \in P\), $y_{0} \in Q$,
		LMO over \(P, Q \subseteq \RR^{n}\)
    \Ensure \(x \in P \cap Q\) or ``disjoint''
		\hrulealg
		\For{$t=0$ \To \dots}
		\State \(u_{t} \gets \argmin_{x \in P} \innp{x_{t} - y_{t}}{x}\)
		\State \(x_{t+1} \gets
		x_{t} + \gamma_{t,1} \cdot (u_{t} - x_{t})\)
		\State \(v_{t} \gets \argmin_{y \in Q} \innp{y_{t} - x_{t+1}}{y}\)
		\State \(y_{t+1} \gets
		y_{t} + \gamma_{t,2} \cdot (v_{t} - y_{t})\)
    \If{\(t = 2^{k}\) for some \(k\)} \label{line:checks}
      \If{$\min_{x \in P} \innp{x_{t+1} - y_{t+1}}{x}
        > \max_{y \in Q} \innp{x_{t+1} - y_{t+1}}{y}$}
        \State \Return ``disjoint'' and
        certificate $\forall x\in P, y\in Q\colon
        \innp{x_{t+1} - y_{t+1}}{x - y} > 0$
      \Else
        \State Solve linear program \eqref{eq:LMO-point}.
        \If{feasible}
          \State \Return a solution $x \in P \cap Q$
        \EndIf
			\EndIf
		\EndIf
		\EndFor
	\end{algorithmic}
\end{algorithm}

\begin{algorithm}
  \caption{Away-step Frank–Wolfe algorithm for intersection}
  \label{alg:away-intersect}
	\begin{algorithmic}[1]
		\Require Points \(x_{0} \in P\), $y_{0} \in Q$,
		LMO over \(P, Q \subseteq \RR^{n}\)
    \Ensure \(x \in P \cap Q\) or ``disjoint''
		\hrulealg
		\For{$t=0$ \To \dots}
		\State \(u_{t} \gets \argmin_{x \in P} \innp{x_{t} - y_{t}}{x}\)
		\State \(v_{t} \gets \argmin_{y \in Q} \innp{y_{t} - x_{t}}{y}\)
    \State Find convex combination
    \((x_{t}, y_{t}) = \sum_{i} \lambda_{i} (w_{i}, z_{i})\)
    with \(w_{i} \in P\) and \(z_{i} \in Q\).
    \label{line:convex-decompose}
    \State \(k \gets \argmax_{i}
    \innp{x_{t} - y_{t}}{w_{i} - z_{i}}\)
    \If{\(\innp{x_{t} - y_{t}}{w_{k} - z_{k} - x_{t} + y_{t}} \geq
      \innp{x_{t} - y_{t}}{x_{t} - y_{t} - u_{t} + v_{t}}\)}
    \State \(a_{t} \gets x_{t} - w_{k}\), \(b_{t} \gets y_{t} - z_{k}\)
    \State \(\gamma_{t} \gets \min\left\{
    \frac{\innp{x_{t} - y_{t}}{w_{k} - z_{k}}}{
      \norm{x_{t} - y_{t}}^{2}}
    - 1,
    \frac{\lambda_{k}}{1 - \lambda_{k}}
    \right\}\)
    \Else
    \State \(a_{t} \gets u_{t} - x_{t}\),
    \(b_{t} \gets v_{t} - y_{t}\)
    \State \(\gamma_{t} \gets \min\left\{
    1 - \frac{\innp{x_{t} - y_{t}}{u_{t} - v_{t}}}{
      \norm{x_{t} - y_{t}}^{2}},
    1\right\}\)
    \EndIf
		\State \(x_{t+1} \gets
    x_{t} + \gamma_{t} \cdot a_{t}\)
		\State \(y_{t+1} \gets
    y_{t} + \gamma_{t} \cdot b_{t}\)
    \If{\(t = 2^{k}\) for some \(k\)} \label{line:away-checks}
      \If{$\min_{x \in P} \innp{x_{t+1} - y_{t+1}}{x}
         > \max_{y \in Q} \innp{x_{t+1} - y_{t+1}}{y}$}
        \State \Return ``disjoint'' and
        certificate $\forall x\in P, y\in Q\colon
        \innp{x_{t+1} - y_{t+1}}{x - y} > 0$
      \Else
        \State Solve linear program \eqref{eq:LMO-point}.
        \If{feasible}
          \State \Return a solution $x \in P \cap Q$
        \EndIf
			\EndIf
		\EndIf
		\EndFor
	\end{algorithmic}
\end{algorithm}

We will now state the guarantees for Algorithm~\ref{alg:dist_LMO_adaptive} for the case of polytopes.
We count solving \eqref{eq:LMO-point} as one block-LMO call.
For comparison, we include guarantees for
the Away-step Frank–Wolfe algorithm \citep[Theorem~1]{lacoste15}
on the product \(P \times Q\)
with the disjointness check of Algorithm~\ref{alg:dist_LMO_adaptive},
i.e., Algorithm~\ref{alg:away-intersect},
a non-alternating variant of
Algorithm~\ref{alg:dist_LMO_adaptive},
where \(x_{t+1}\) and \(y_{t+1}\) are chosen
with the help of
\(\argmin_{x \in P} \innp{x_{t} - y_{t}}{x}\) and
\(\argmin_{y \in Q} \innp{y_{t} - x_{t}}{y}\),
i.e., two block-LMO calls.

\begin{proposition}[Adaptive variant]
  \label{prop:adaptive}
  Let $P, Q$ be polytopes with diameters \(D_{p}\) and \(D_{Q}\),
  respectively.
  Algorithm~\ref{alg:dist_LMO_adaptive} has the following guarantees,
  with an additional number of LMO calls logarithmic in the number of
  block-LMO calls:
\begin{enumerate}
	\item If $P \cap Q \neq \emptyset$, then after no more than 
\[
\frac{16 (1 + 2 \sqrt{2}) (D_{P}^{2} + D_{Q}^{2})}{\varepsilon_{PQ}^2}
\]
block-LMO calls,
the algorithm returns $x \in P \cap Q$ with the agnostic step-size
rule and no more than
\[
\frac{32 ((D_{P} + D_{Q}) \cdot \max \{D_{P}, D_{Q}\}
  + 2 (D_{P}^{2} + D_{Q}^{2}))
}{\varepsilon_{PQ}^2}
\]
block-LMO calls with the short step-size rule.
\item If $P \cap Q = \emptyset$, then after no more than
  \[
    \frac{27 (1 + 2 \sqrt{2})
      (D_{P}^{2} + D_{Q}^{2})}{\dist(P,Q)^{2}}
\]
block-LMO calls the algorithm certifies $P \cap Q = \emptyset$
with the agnostic step-size rule and no more than
\[
  \frac{64}{\dist(P,Q)^{2}}
    \left(
      (D_{P} + D_{Q} + \dist(P,Q)) \cdot \max \{D_{P}, D_{Q}\}
      + 2 (D_{P}^{2} + D_{Q}^{2})
    \right)
 \]
block-LMO calls with the short step-size rule.
\end{enumerate}
\begin{proof}
If $P \cap Q \neq \emptyset$ and using the agnostic step size rule,
then by Proposition~\ref{prop:setIntersection},
after \(T \doteq 4 (1 + 2 \sqrt{2}) (D_{P}^{2} + D_{Q}^{2})
/ \varepsilon_{PQ}^{2} - 2\) iterations
one has \(\norm{x_{t} - y_{t}} < \varepsilon_{PQ}\)
for all \(t \geq T\),
and hence by Remark~\ref{rem:exact} solving \eqref{eq:LMO-point}
provides a point \(x\) in the intersection.
Due to delayed tests, the algorithm executes at most \(2 T\)
iterations, with two block-LMO calls per iteration as before.

The other bounds are obtained similarly.
The case $P \cap Q = \emptyset$
follows via Corollary~\ref{cor:disjoint-par-free},
where the bound on block-LMO calls is twice the bound on iterations.
\end{proof}
\end{proposition}

The proofs of Proposition~\ref{prop:setIntersection}, Corollary~\ref{cor:disjoint}, Proposition~\ref{prop:adaptive} and Corollary~\ref{cor:disjoint-par-free} heavily rely on the Cyclic Block-Coordinate Conditional Gradient method from \citet{beck2015cyclic} (see Theorem~\ref{thm:almo}) and its convergence rate. In the corresponding proof one has to account for errors arising from using the gradient on a single component only rather than the gradient on the whole product, this error accumulation in the current proof limits the overall convergence rate to $\mathcal O(1/t)$. If one is willing to give up on alternations, solve a single linear program over the product of the two sets in each iteration, and use e.g., the Away-step Frank-Wolfe algorithm \citep{wolfe70}, one vastly improve the convergence rate in the polytope setting as we show in the following remark. Also any of the other variant achieving linear convergence for polytopes is fine, such as e.g., Fully-corrective FW, FW with local linear minimization \citep{garber2016linearly}, Blended Conditional Gradients \citep{pok18bcg}, or Blended Pairwise Conditional Gradients \citep{TTP2021}; see \citet{CGFWSurvey2022} for an overview; we formulate the variant with Away-step Frank-Wolfe algorithm below. We will use that the pyramidal width of \(P \times Q\) in the Euclidean norm is
\(\sqrt{\delta_{P}^{2} \delta_{Q}^{2}
  / (\delta_{P}^{2} + \delta_{Q}^{2})}\)
by \cite[Theorem~4.5]{IWCRP2025},
where \(\delta_{P}\) and \(\delta_{Q}\) are the pyramidal widths of
\(P\) and \(Q\) in the Euclidean norm, respectively.

\begin{remark}[Improved rates without alternations]
  \label{rem:adaptive_afw}
  Let $P, Q$ be polytopes with diameters \(D_{p}\), \(D_{Q}\) and
  pyramidal widths \(\delta_{P}\), \(\delta_{Q}\), respectively.
  Then Algorithm~\ref{alg:away-intersect}
  has the following guarantees:
\begin{enumerate}
\item If $P \cap Q \neq \emptyset$, then
  it returns a point in the intersection
  after no more than
\begin{equation}
  \label{eq:3}
  6 +
  24
  \frac{(\delta_{P}^{2} + \delta_{Q}^{2}) (D_{P}^{2} + D_{Q}^{2})}
  {\delta_{P}^{2} \delta_{Q}^{2}}
  \ln \frac{D_{P}^{2} + D_{Q}^{2}}{\varepsilon_{P, Q}^{2}}
\end{equation}
block-LMO calls.
\item If $P \cap Q = \emptyset$, then
  it certifies emptiness of intersection
  after no more than
\begin{equation}
  \label{eq:4}
    6 +
    48
    \frac{(\delta_{P}^{2} + \delta_{Q}^{2}) (D_{P}^{2} + D_{Q}^{2})}
    {\delta_{P}^{2} \delta_{Q}^{2}}
    \ln \frac{\dist(P, Q)}{2}
\end{equation}
block-LMO calls.
\end{enumerate}
\end{remark}

In the case of polytopes, we suspect that one can improve the convergence rates of schemes such as the Cyclic Block-Coordinate Conditional Gradient method when paired with steps of linearly convergence variants. We leave this as an open problem.

\section{Preliminary computation}
\label{sec:experiments}

The computation was done
on a machine with 512GB or 1024GB RAM
and 48-core Intel(R) Xeon(R) Gold 6342 CPU
operating at frequency 2.80 GHz.
We have used the \texttt{FrankWolfe.jl} Julia package
(version 0.6.3) with \texttt{Julia} (version 1.12.5)
on the operating system Debian GNU/Linux (version 12.13). The following computational experiments should be considered as indicative while a full numerical benchmark is beyond the scope of this work; we refer the reader to the \texttt{FrankWolfe.jl} package\footnote{see \url{https://zib-iol.github.io/FrankWolfe.jl/stable/basics/\#Alternating-Linear-Minimization-(ALM)}} for comprehensive benchmarking \citep{BDHKPTVW2024}.

We ran several representative experiments on various configurations and choices of feasible sets to compare ALM against POCS. The arising projection problems in POCS, in and of themselves non-trivial quadratic optimization problems, were solved with the same FW variants as the one used for ALM. We considered two common variants of the Frank-Wolfe algorithm to implement ALM (and the projections in POCS): the standard Frank-Wolfe algorithm with lazification (LFW) (without lazification a single POCS iteration might take roughly as long as the full ALM run; see \citet{pok17lazy}) as well as the Blended Pairwise Conditional Gradient algorithm (BPCG) (the current default in the \texttt{FrankWolfe.jl} package; see \citet{TTP2021}). We ran both algorithms, ALM and POCS until the Frank-Wolfe gap was less than $10^{-7}$, the precision for the projection subproblem was set to $10^{-8}$ to ensure that the subproblems are solved to sufficient accuracy in line with conditional gradient sliding approaches; see \citep{lan2016conditional}, see also \citep{CGFWSurvey2022}. In more elaborate implementations also an adaptive error criterion for the projection subproblems could be used, which however is beyond the scope. Note that the Frank-Wolfe gap can also be evaluated as a stopping criterion for POCS. An overview of the experiments can be found in Table~\ref{tab:experiments}. We used the Birkhoff polytope as the second set for all experiments, given its prevalence as a popular example of a polytope with non-trivial geometry in the Frank-Wolfe literature; see \citet{CP2021} for a discussion on the complexity of linear minimization vs projections.

As can be seen in many cases ALM outperforms POCS both in time and
iterations, while in others POCS is faster. Note, that sometimes it
looks as if POCS is running beyond the accuracy level, which is simply
because the algorithm has made a projection before verifying the
stopping criterion. We would like to point out that POCS heavily
benefits from the fast FW variants employed for the subsolver as
already observed in, e.g., \citep{HRBDP2025,LHPD2025}, which exploit
BPCG as subsolver for projections. If we would have used standard
Frank-Wolfe without lazification for the projection subproblems, POCS
would have been several orders of magnitude slower.

\begin{table}[ht]
  \centering
  \caption{Overview of computational experiments. All instances use the Birkhoff polytope for $Q$. The dimension is $n=100$.}
  \label{tab:experiments}
  \begin{tabular}{@{}llcc@{}}
    \toprule
    Set $P$ & Set $Q$ & Disjoint & Intersecting \\
    \midrule
    $\ell_2$-ball        & Birkhoff polytope & \cref{fig:l2ball-birkhoff-09} & \cref{fig:l2ball-birkhoff-11} \\
    Nuclear norm ball    & Birkhoff polytope & \cref{fig:nucnorm-birkhoff-09} & \cref{fig:nucnorm-birkhoff-11} \\
    Spectrahedron        & Birkhoff polytope & --- & \cref{fig:spectrahedron-birkhoff-11} \\
    \bottomrule
  \end{tabular}
\end{table}

\begin{figure}
  \centering
  \includegraphics[width=.49\textwidth]{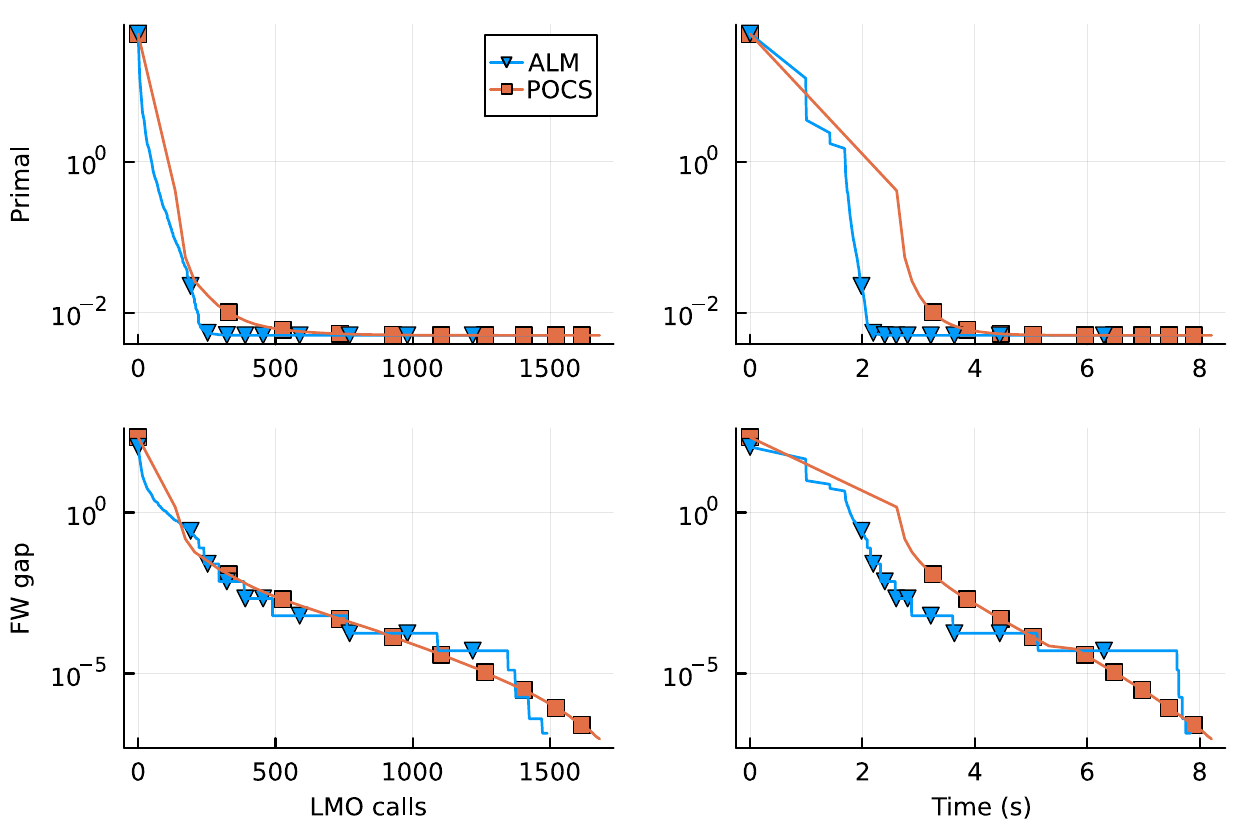}%
  \hfill
  \includegraphics[width=.49\textwidth]{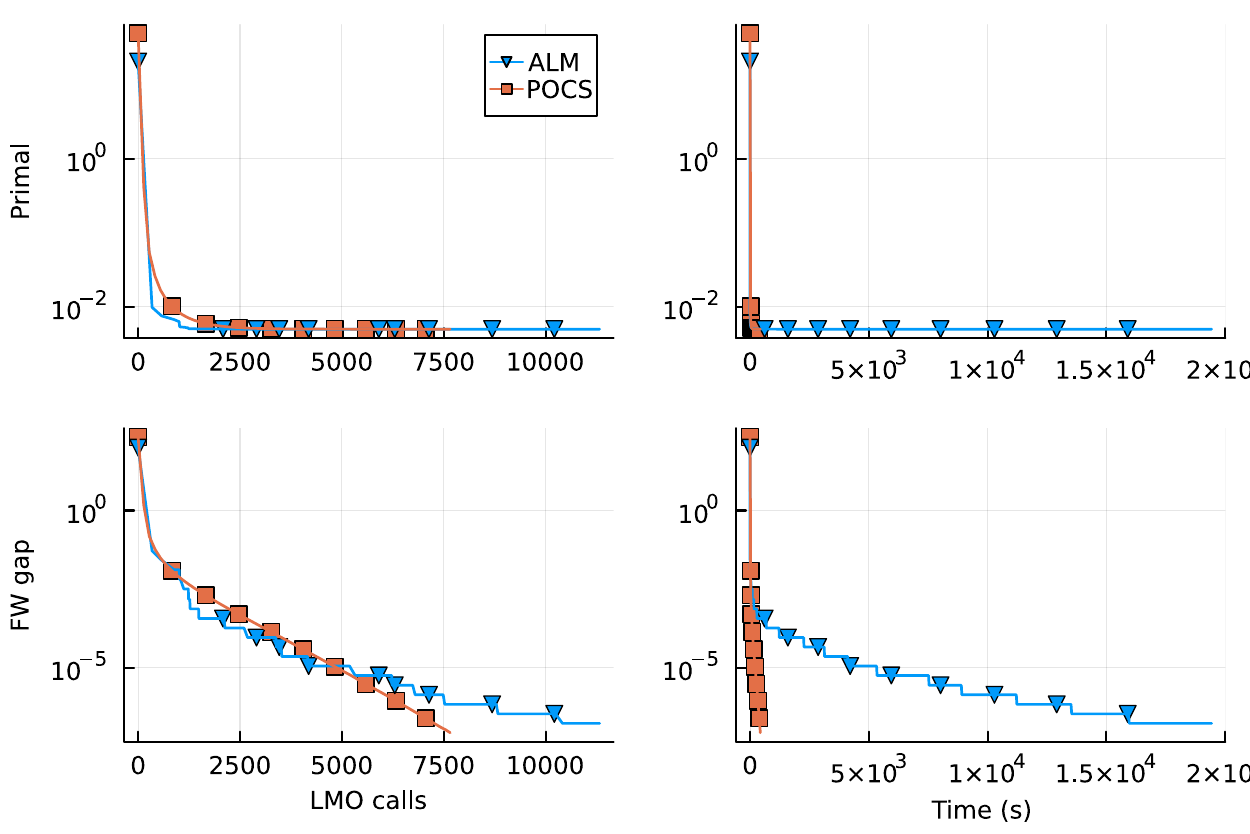}
  \caption{$\ell_2$-ball vs.\ Birkhoff polytope, disjoint. Left: BPCG. Right: LFW.}
  \label{fig:l2ball-birkhoff-09}
\end{figure}
\begin{figure}
  \centering
  \includegraphics[width=.49\textwidth]{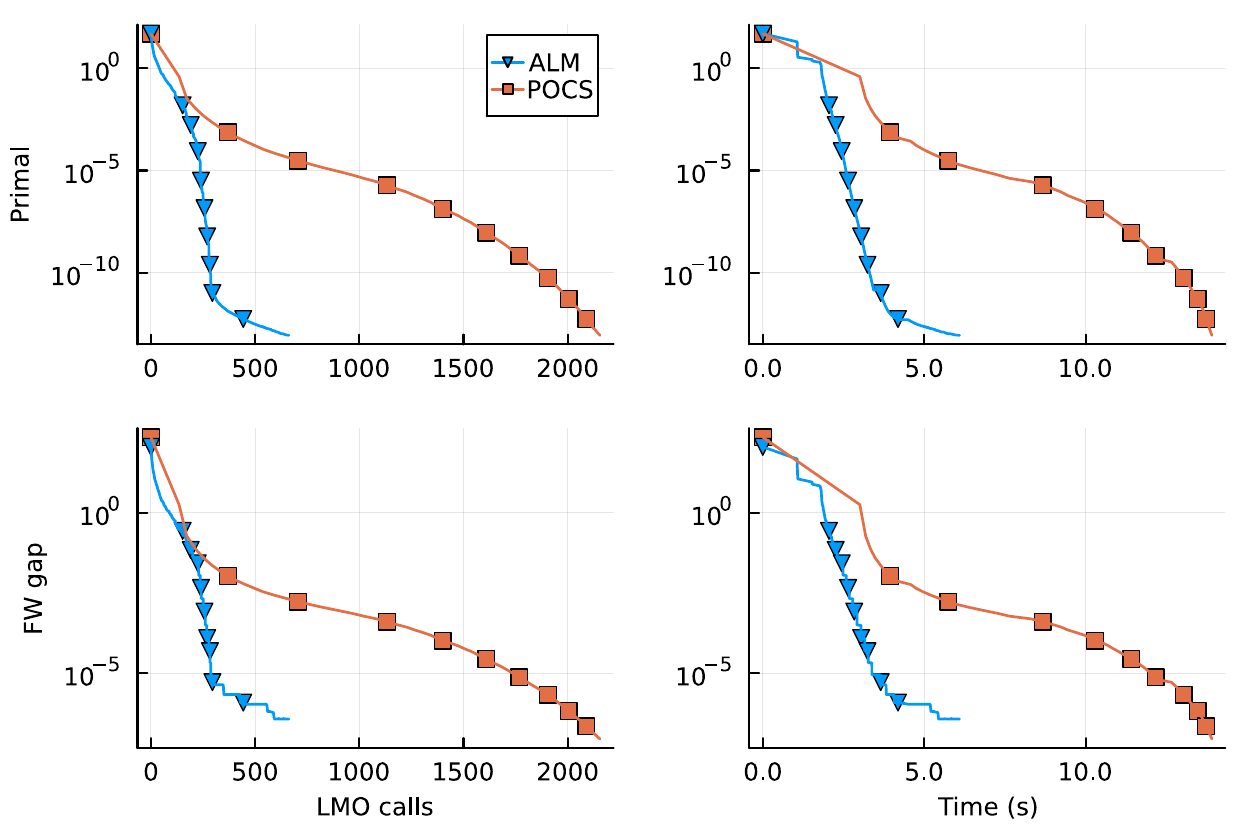}%
  \hfill
  \includegraphics[width=.49\textwidth]{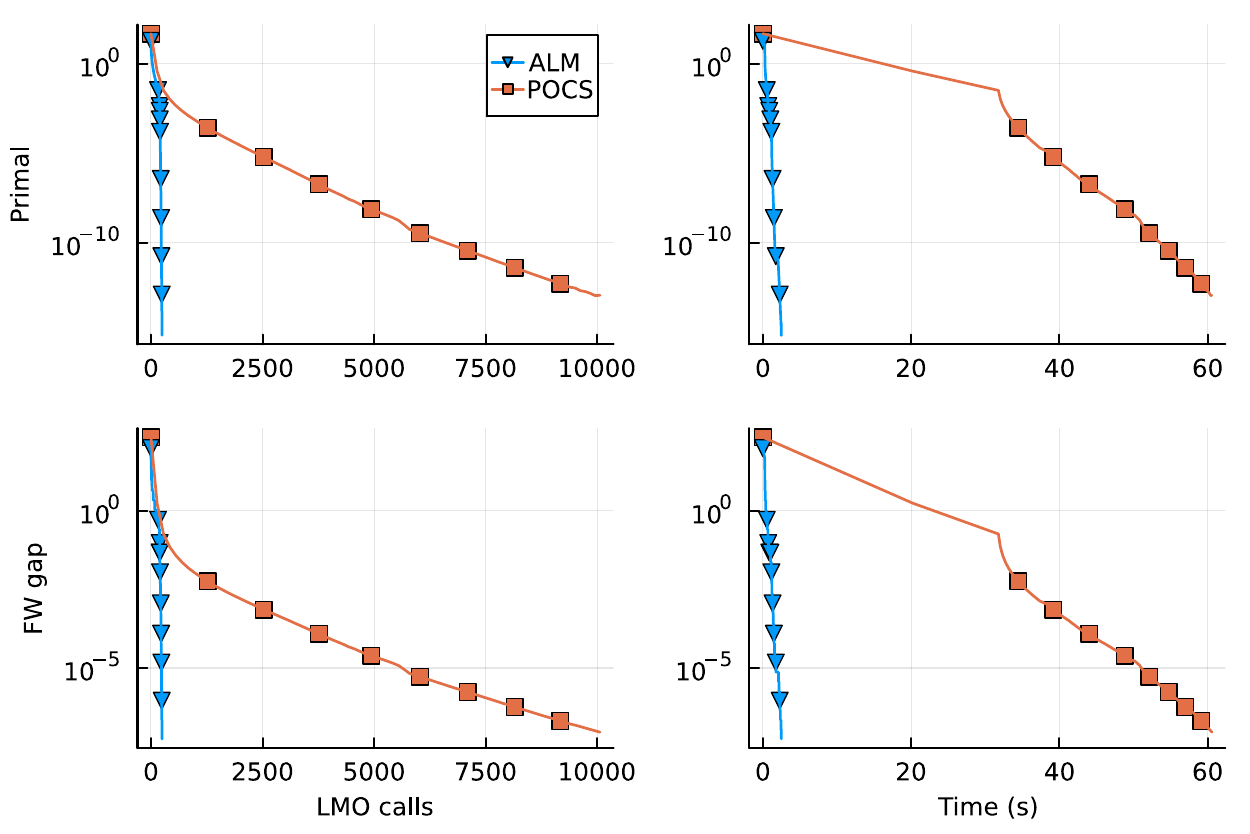}
  \caption{$\ell_2$-ball vs.\ Birkhoff polytope, intersecting. Left: BPCG. Right: LFW.}
  \label{fig:l2ball-birkhoff-11}
\end{figure}
\begin{figure}
  \centering
  \includegraphics[width=.49\textwidth]{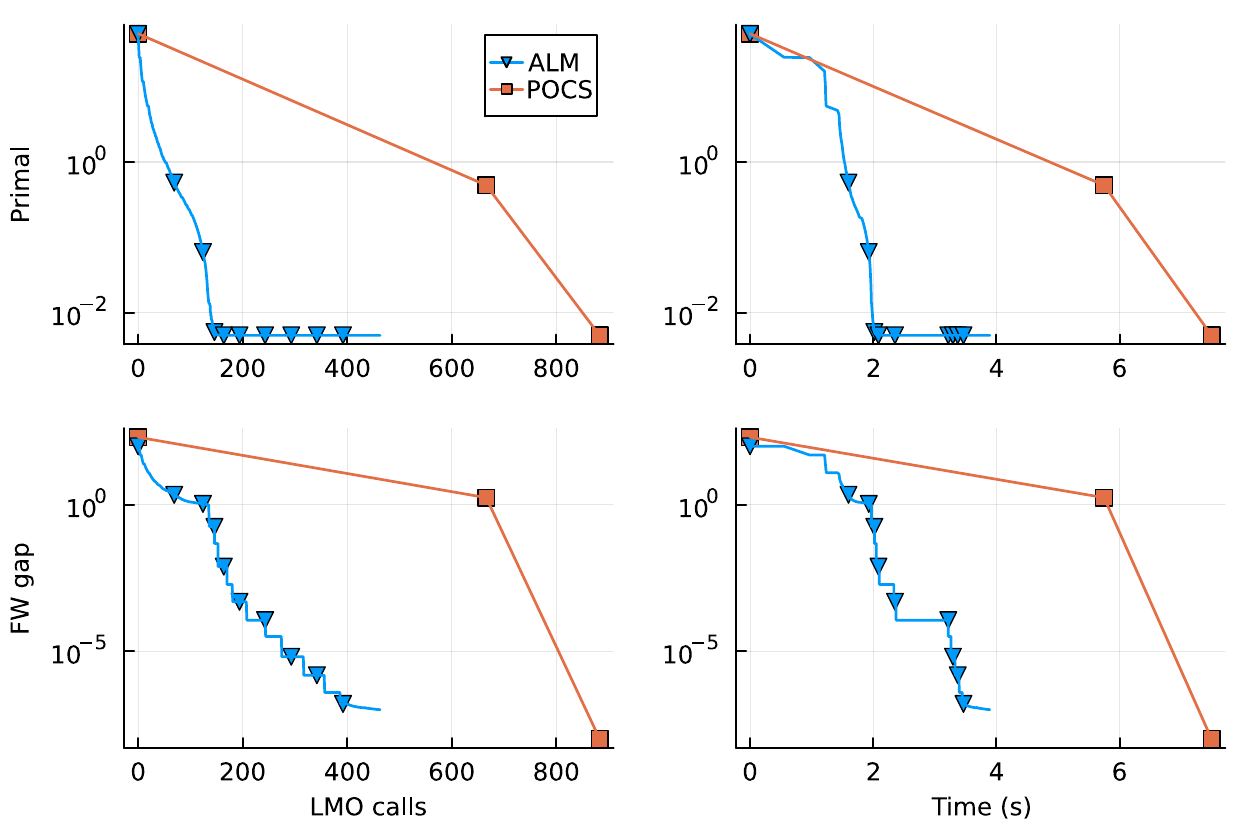}%
  \hfill
  \includegraphics[width=.49\textwidth]{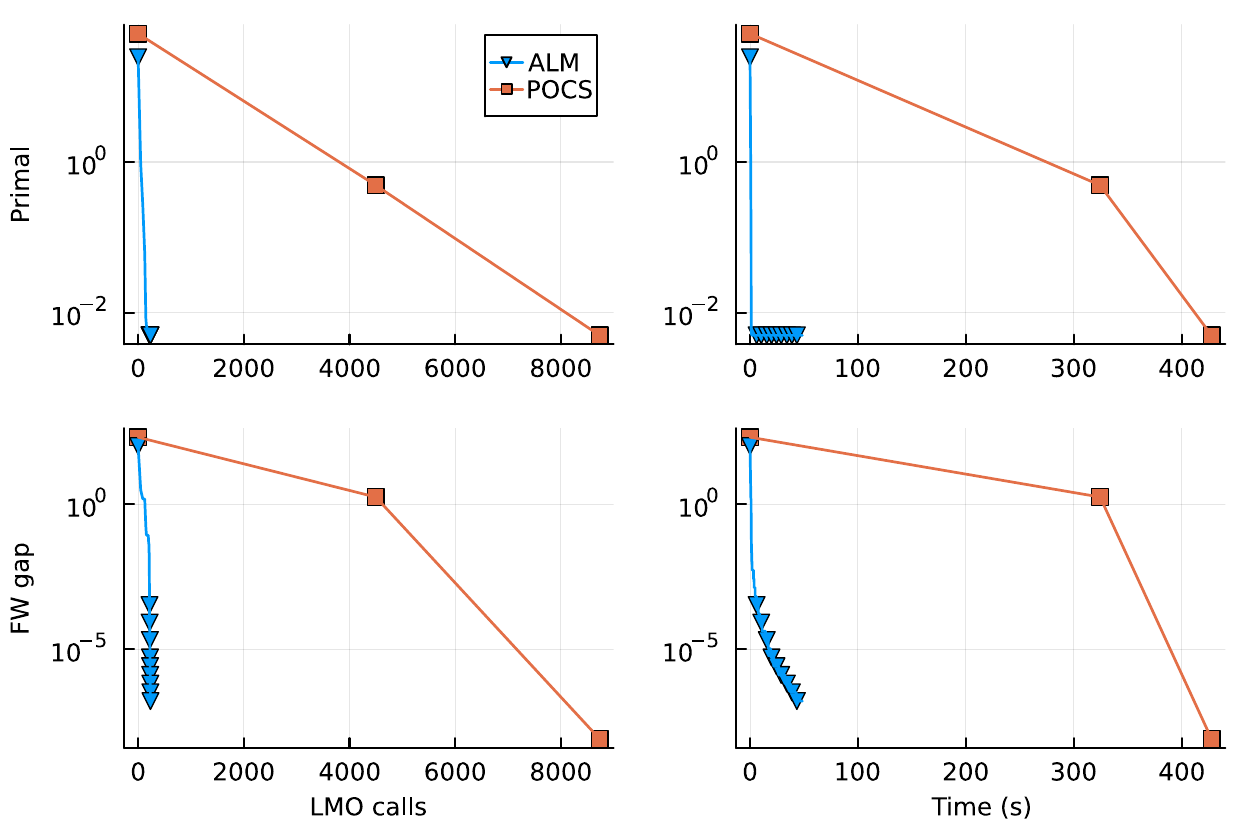}
  \caption{Nuclear norm ball vs.\ Birkhoff polytope, disjoint. Left: BPCG. Right: LFW.}
  \label{fig:nucnorm-birkhoff-09}
\end{figure}
\begin{figure}
  \centering
  \includegraphics[width=.49\textwidth]{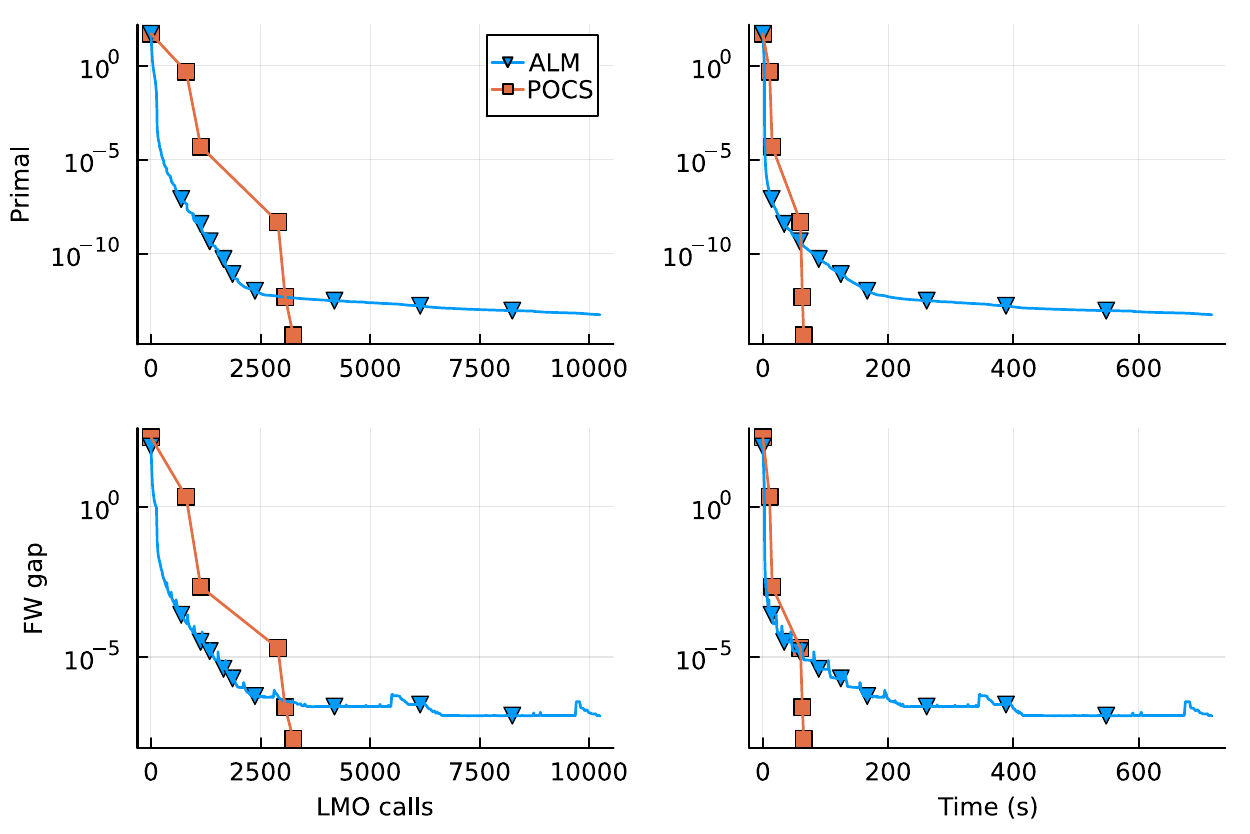}%
  \hfill
  \includegraphics[width=.49\textwidth]{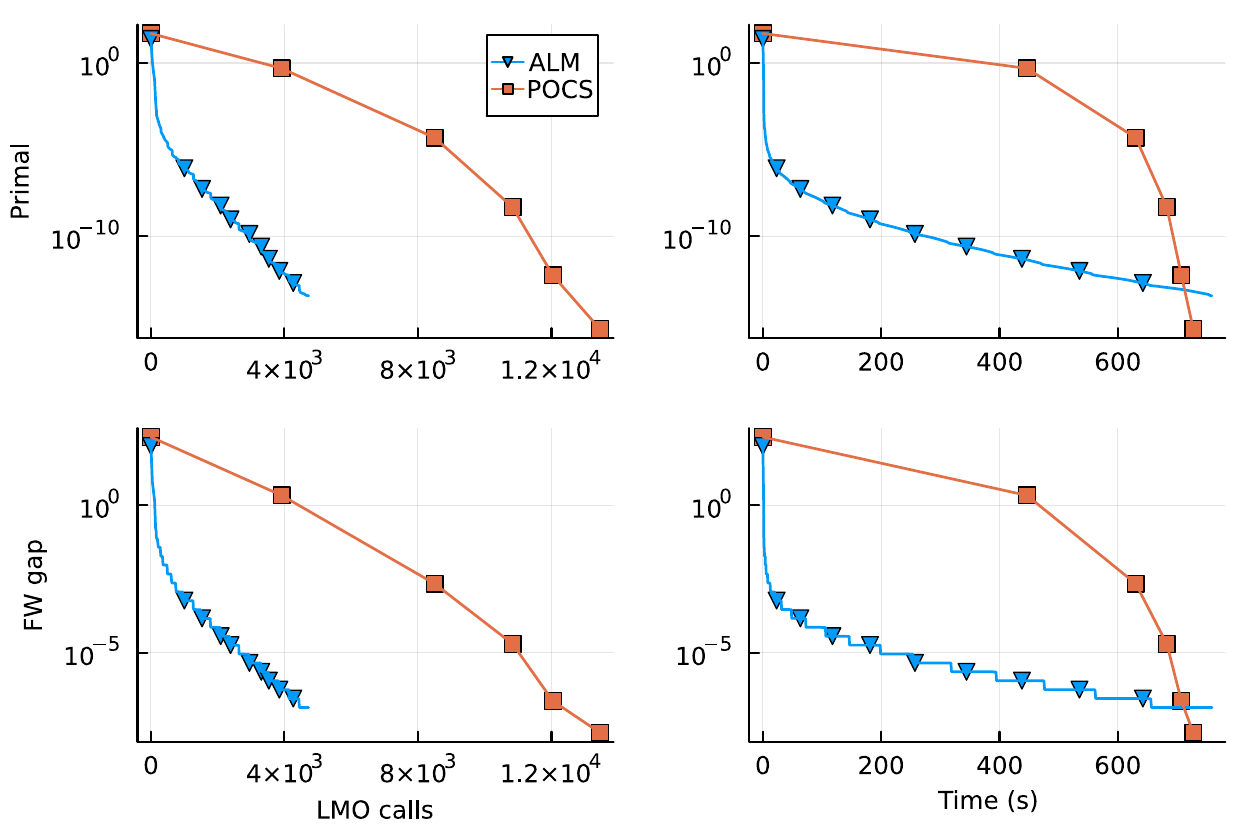}
  \caption{Nuclear norm ball vs.\ Birkhoff polytope, intersecting. Left: BPCG. Right: LFW.}
  \label{fig:nucnorm-birkhoff-11}
\end{figure}
\begin{figure}
  \centering
  \includegraphics[width=.49\textwidth]{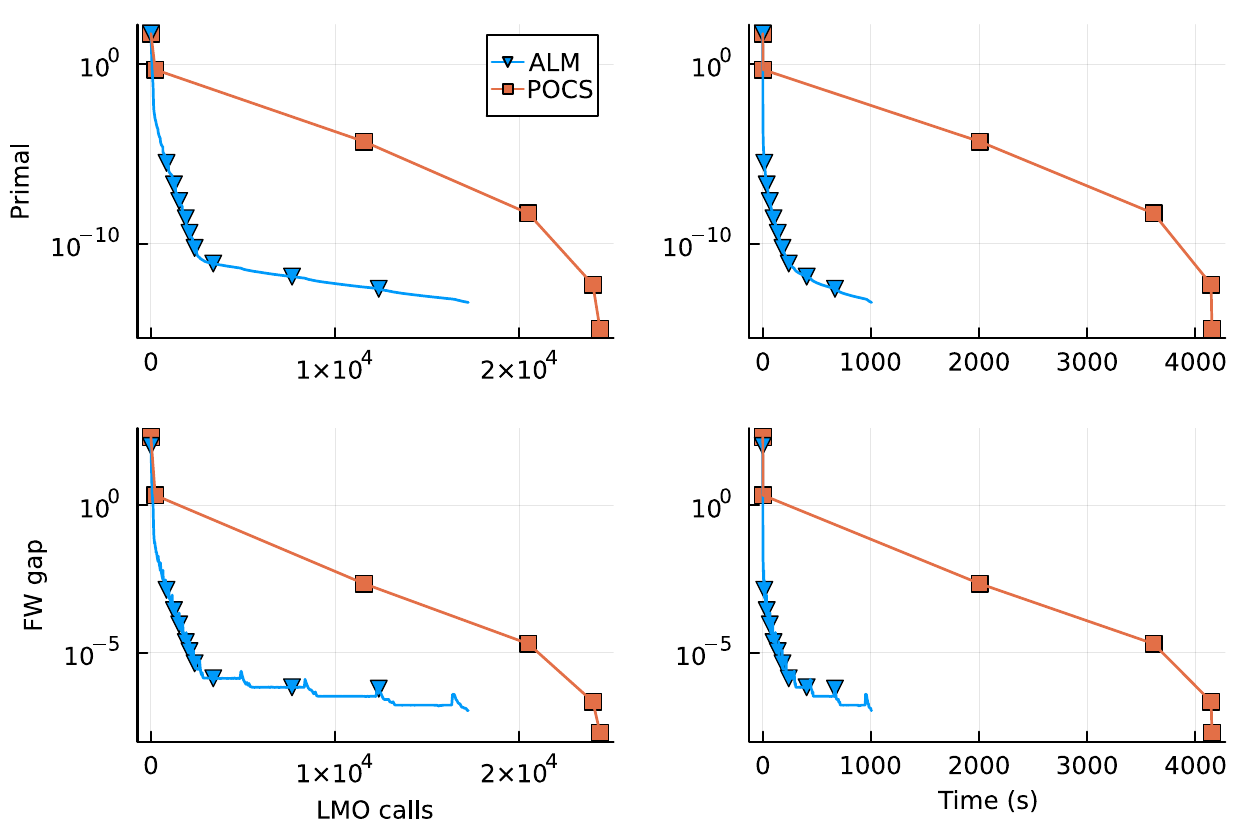}%
  \hfill
  \includegraphics[width=.49\textwidth]{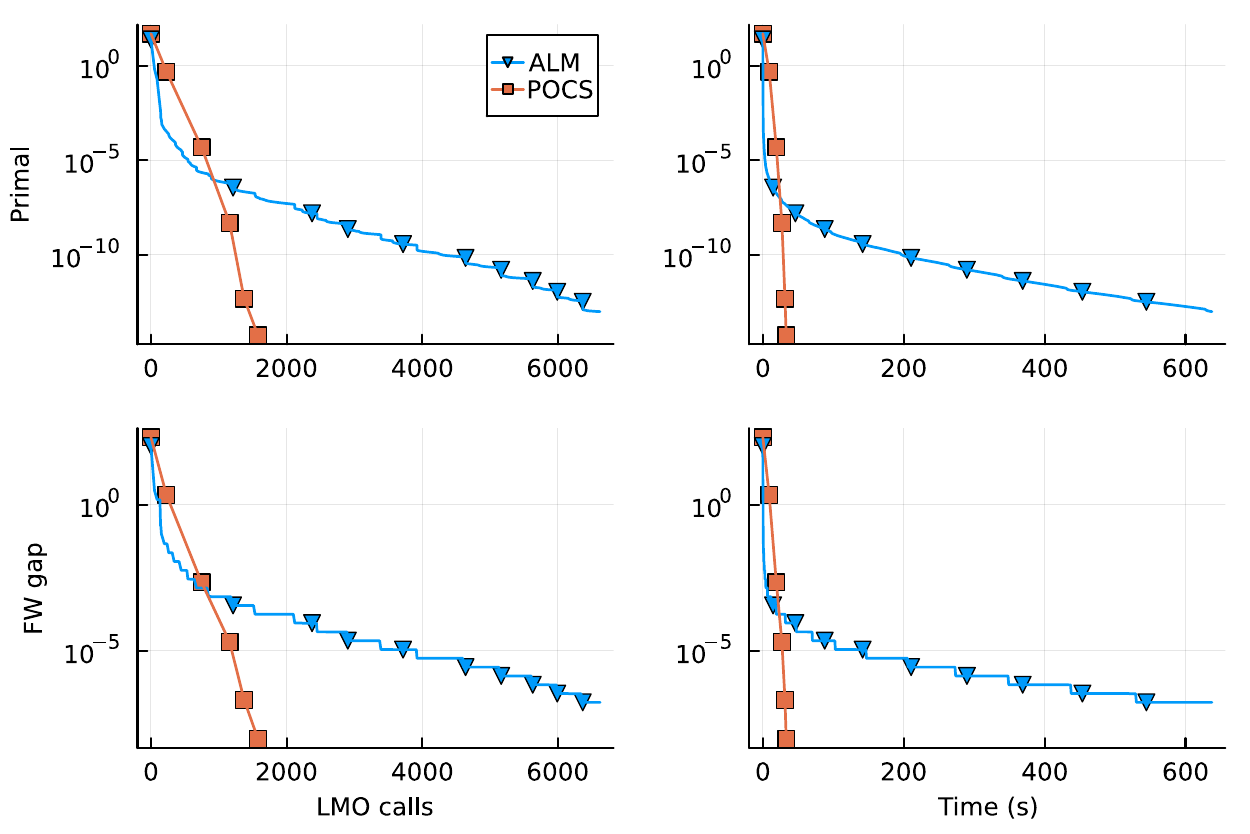}
  \caption{Spectrahedron vs.\ Birkhoff polytope, intersecting. Left: BPCG. Right: LFW.}
  \label{fig:spectrahedron-birkhoff-11}
\end{figure}
\clearpage

\section*{Acknowledgements}

Research reported in this paper was partially supported by the DFG Cluster of Excellence MATH+ (EXC-2046/1 and EXC-2046/2, project id 390685689) funded by the Deutsche Forschungsgemeinschaft (DFG). The third author was supported by the Einstein Foundation Berlin.


\bibliography{refs}
\bibliographystyle{icml2021}

\end{document}